\documentclass[fontsize=12pt,a4paper,headings=normal,
twoside=false,leqno,parskip=half-,abstract=true]{scrartcl}
\usepackage[english]{babel}
\usepackage[utf8]{inputenc}
\usepackage{hyperref}
\hypersetup{
 pdftitle={SymmSphere},
 pdfauthor={Phillipo Lappicy},
 colorlinks=true,
 linkcolor=blue,
 citecolor=blue,
 filecolor=blue,
 urlcolor=blue}

\usepackage{graphicx}
\usepackage[format=plain,labelfont=bf,font=small]{caption}
\usepackage{xcolor}
\usepackage[arrow, matrix, curve]{xy}
\usepackage{float}

\usepackage{caption}
\usepackage{tabulary}
\usepackage{array}
\newcolumntype{N}[1]{>{\centering\arraybackslash}m{#1}}

\usepackage{amsmath,amsthm}
\usepackage{amssymb} 

\makeatletter
\newcommand{\tpitchfork}{%
  \vbox{
    \baselineskip\z@skip
    \lineskip-.52ex
    \lineskiplimit\maxdimen
    \m@th
    \ialign{##\crcr\hidewidth\smash{$-$}\hidewidth\crcr$\pitchfork$\crcr}
  }%
}
\makeatother
\usepackage{latexsym}

\usepackage[notref,notcite,color,final 
]{showkeys}

\definecolor{refkey}{rgb}{1,0,0}
\definecolor{labelkey}{rgb}{1,0,0}

\usepackage{tikz}

\usepackage[textwidth=2cm,textsize=small,backgroundcolor=none]{todonotes}

  \mathchardef\ordinarycolon\mathcode`\:
  \mathcode`\:=\string"8000
  \begingroup \catcode`\:=\active
    \gdef:{\mathrel{\mathop\ordinarycolon}}
  \endgroup

\newtheorem{thm}{Theorem}[section]
\newtheorem{lem}[thm]{Lemma}

\newenvironment{pf}[1][Proof]{\begin{trivlist}
\item[\hskip \labelsep {\bfseries #1}]}{\end{trivlist}}

\usepackage{comment}
\usepackage{caption}

\begin{document}

\title{{\LARGE{A symmetry property for elliptic equations on the sphere}}}

\author{
 \\
{~}\\
Phillipo Lappicy*\\
\vspace{2cm}}

\maketitle
  \thispagestyle{empty}

\vfill

$\ast$\\
Instituto de Ciências Matemáticas e de Computação\\
Universidade de S\~ao Paulo\\
Avenida trabalhador são-carlense 400\\
13566-590, São Carlos, SP, Brazil\\


\newpage
\pagestyle{plain}
\pagenumbering{arabic}
\setcounter{page}{1}

\begin{abstract}

The goal of this paper is to study how the symmetry of the spherical domain influences solutions of elliptic equations on such domain. The method pursued is a variant of the moving plane method, discovered by Alexandrov (1962) and used for differential equations by Gidas, Ni and Nirenberg (1979). We obtain a reflectional symmetry result with respect to maxima and minima of solutions.

\ 


\end{abstract}


\section{Main result}\label{cpt:sym}

\numberwithin{equation}{section}
\numberwithin{figure}{section}
\numberwithin{table}{section}

Consider the fully nonlinear elliptic equation
\begin{equation}\label{symPDE}
    0 = F(u,\Delta_{\mathbb{S}^2}u)
\end{equation}
where the spatial domain is $(\theta,\phi) \in {\mathbb{S}^2}$, $\Delta_{\mathbb{S}^2}$ is the Laplace-Beltrami on the sphere. Assume that $F$ is analytic, and due to \cite{CaoRammahaTiti00}, $u$ is also analytic. Assume also that $F$ is uniformly elliptic, that is, $F_q\geq \delta >0$ everywhere, with $q:=\Delta_{\mathbb{S}^2}u$ and some $\delta\in\mathbb{R}_+$.



We seek a symmetry property for solutions $u$ of the equation \eqref{symPDE}, which is different than the celebrated Gidas, Ni and Nirenberg \cite{GidasNiNirenberg}, but still uses the well known moving plane method developed by Alexandrov \cite{Alexandrov62} and further by Serrin \cite{Serrin71}. 

The Gidas, Ni and Nirenberg result displays a symmetry property of non-negative solutions. We, on the other hand, desire to characterize sign-changing solutions. For example, most solutions change sign for the elliptic Chafee-Infante problem given by $F(u,q)=q+\lambda u[1-u^2]$, with $\lambda\in\mathbb{R}$ and Dirichlet boundary conditions. Indeed, the only solution that does not change sign is the trivial solution $u\equiv 0$ which is clearly radially symmetric. On the other hand, there are several bifurcating solutions from $u\equiv 0$ that change sign, and their properties should be developed.

The function $u(\theta,\phi)$ has \emph{axial extrema} if its extrema occur as axis from the north to south pole. Mathematically, if $u_\phi(\theta_0,\phi_0)=0$ for a fixed $(\theta_0,\phi_0)\in \mathbb{S}^2$, then $u_\phi(\theta,\phi_0)=0$ for any $\theta\in[0,\pi]$. In that case, the extrema depend only at the position in $\phi$. Note that there are finitely many axial extrema, since $u$ is analytic. Denote them by $\{ \phi_i \}_{i=1}^N$.

Axial extrema are \emph{leveled} if all axial maxima $\phi_j$ attain the same value $u(\theta,\phi_j)=M(\theta)$, and all axial minima $\phi_k$ attain the same value $u(\theta,\phi_k)=m(\theta)$.

As a first step, we consider solutions with leveled axial extrema, which are not within the scope of the usual symmetry property for non-negative solutions. The symmetry obtained is of reflectional type, where the reflection of $\phi$ with respect to the extrema $\phi_i$ is described by 
\begin{equation}\label{Rdef}
    R_{\phi_i}(\phi):=2\phi_i-\phi.
\end{equation}  

We mention that, the bifurcating solutions of the Chafee-Infante problem can be computed numerically, as in \cite{UeckerWetzelRademacher14} and references therein. In those cases, there are some examples of solutions which have \emph{leveled axial extrema}.

The main result is described below.

\begin{thm}\emph{\textbf{Reflectional symmetry}} \label{symhetthm4}

Suppose that $u(\theta,\phi)$ is a non constant solution of \eqref{symPDE} such that all its extrema are leveled and axial at $\{ \phi_i \}_{i=1}^N$. Then $\phi_i=(\phi_{i-1}+\phi_{i+1})/2$ and
\begin{equation}\label{SYM3}
    u(\theta,\phi)=u(\theta,R_{\phi_i}(\phi))
\end{equation}
for all $i=1,...,N$, where $(\theta,\phi)\in[0,\pi] \times [\phi_{i-1},\phi_{i}]$, $\phi_{0}:=\phi_{N}$ and $\phi_{N+1}:=\phi_1$.
\end{thm}

Positive solutions of elliptic equations on a ball with Dirichlet boundary conditions have a symmetry property obtained by Gidas, Ni and Nirenberg \cite{GidasNiNirenberg}. Proving a similar symmetry property for positive solutions in the sphere has some difficulties. In particular, the domain has no boundary and it is not clear where to start the moving plane method. 

This problem was partially solved by Padilla \cite{Padilla} for particular convex subsets of the sphere, and later considered by Kumaresan and Prajapat \cite{KumaresanPrajapat98} for subsets of the sphere contained in a hemisphere. Later, Brock and Prajapat \cite{BrockPrajapat00} considered subsets containing hemisphere, but still not the full sphere. All these methods rely on a stereographic projection, so that domains within the sphere are transformed into domains in the Euclidean space and one can apply the moving plane method. 

Hence, the usual symmetry property in the sphere is still an open problem: are non-negative solutions of an elliptic problem on the sphere axially symmetric, depending only on the angle $\theta$?

One could try such a stereographic projection in order to apply one of Frankel's \cite{Fraenkel00} collection of symmetry results for the whole plane $\mathbb{R}^2$. In this case, it is required admissible decay conditions outside large sets in the plane in order to start the moving plane method. These conditions could possibly be satisfied due to the singularity at the north pole, which imposes Neumann boundary conditions, and hence growth conditions outside large sets. Another possibility to obtain axisymmetry of solutions is to adapt the work of Pacella and Weth \cite{PacellaWeth07} or Salda{\~{n}}a and Weth \cite{SaldanaWeth12}. 

On another note, we mention that the proof of Theorem \ref{symhetthm4} can be replicated for the problem in the ball with \emph{any} boundary condition, yielding a symmetry result for sign changing solutions. Therefore, Theorem \ref{symhetthm4} displays symmetries of sign-changing solutions which are not covered in the usual Gidas, Ni and Nirenberg case. Nevertheless, we focused on the problem in the sphere, since this equation models the event horizon of certain black holes, as in \cite{LappicyBlackHoles}. 

\begin{figure}[ht]\centering
\begin{tikzpicture}[scale=1]
    \draw[->] (-2.2,0) -- (6.5,0) node[right] {$\phi$};
    \draw[->] (-2,-1.1) -- (-2,1.1) node[above] {$u(\theta,\phi)$};
    
    \filldraw [black] (-2,0) circle (1pt) node[anchor=north east]{$0$};
    \filldraw [black] (-2,1) circle (1pt) node[anchor= east]{$M(\theta)$};
    \filldraw [black] (-2,-1) circle (1pt) node[anchor= east]{$m(\theta)$};

    \draw (-2,0) sin (-1.5,1) cos (-1,0);
    \draw (-1,0) sin (-0.5,-1) cos (0,0);
    \filldraw [black] (0,0) circle (1pt) node[anchor=north west]{$\frac{\pi}{2}$};
    
    \draw (0,0) sin (0.25,1) cos (0.5,0);
    \draw (0.5,0) sin (0.75,-1) cos (1,0);
    \draw (1,0) sin (1.25,1) cos (1.5,0);
    \draw (1.5,0) sin (1.75,-1) cos (2,0);
    \filldraw [black] (2,0) circle (1pt) node[anchor=north west]{$\pi$};   
    
    \draw (2,0) sin (2.5,1) cos (3,0);
    \draw (3,0) sin (3.5,-1) cos (4,0);
    \filldraw [black] (4,0) circle (1pt) node[anchor=north west]{$\frac{3\pi}{2}$};  
    
    \draw (4,0) sin (4.25,1) cos (4.5,0);
    \draw (4.5,0) sin (4.75,-1) cos (5,0);
    \draw (5,0) sin (5.25,1) cos (5.5,0);
    \draw (5.5,0) sin (5.75,-1) cos (6,0);
    \filldraw [black] (6,0) circle (1pt) node[anchor=north west]{$2\pi$};   
   
\end{tikzpicture}
\caption{Function with biradial symmetry $C_{2v}$}
\label{cyclicsym}
\end{figure}
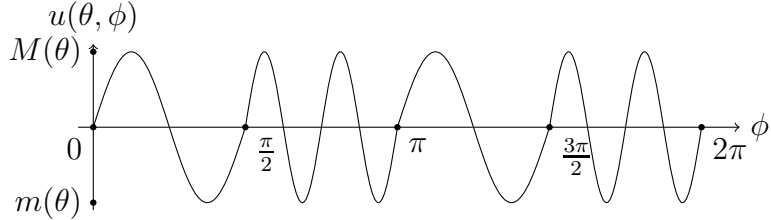

 Lastly, we mention that the hypothesis of the main result do not imply dihedral symmetry. Indeed, one can consider the cyclic symmetric group $C_{2v}=\mathbb{Z}_2\times \mathbb{Z}_2$, in Schönflies notation. Consider a sphere with 4 axis (from the north to south pole) splitting the sphere in 4 parts. Such group acts by reflecting two opposite quadrants. In biology this is called biradial symmetry, and Ctenophora are living beings exhibiting such symmetry. We depict a $\pi$-periodic function that exhibit such symmetry in Figure \ref{cyclicsym}. Therefore, Theorem \ref{symhetthm4} displays more symmetry that is not apparent in the equivariant branching lemma developed by Vanderbauwhede \cite{Vanderbauwhede82}, which shows that the bifurcating branch has the same symmetry as the eigendirection such bifurcation is occuring.


The remaining is organized as follows. In Section \ref{sec:moving}, we provide a proof of Theorem \ref{symhetthm4}, and Section \ref{sec:conclusion} we discuss possible generalizations of Theorem \ref{symhetthm4}.

\section{Proof}\label{sec:moving}

The starting point of the moving arc method is at some axial extrema $\phi_{i-1}$. Without loss of generality, we can assume $\phi_{i-1}=0$. Otherwise, consider $\phi\mapsto \phi-\phi_{i-1}$. 

Define a small sector within the sphere nearby $\phi_{i-1}=0$,
\begin{equation*}
    \Omega_\epsilon :=\{(\theta,\phi)\in \mathbb{S}^2 \text{ $|$ } \theta\in ,[0,\pi]\phi\in (0,\epsilon) \}
\end{equation*}
with small $\epsilon>0$. The boundary $\partial\Omega_\epsilon$ is given by two axis, namely at $\phi=0$ and $\phi=\epsilon$. We call the latter by \emph{$\epsilon$-arc}. We also consider, for instance, $\Omega_{-\epsilon}$ a sector in the other direction of $\phi$.

Then, we will move the $\epsilon$-arc by increasing $\epsilon$, and consider reflections along it, defined by $R_\epsilon (\phi):=2\epsilon-\phi$. Due to the  periodic boundary conditions in $\phi$, 
one has $\Omega_\epsilon,R_{\epsilon}(\Omega_\epsilon)\subset \mathbb{S}^2$ for arbitrarily large $\epsilon$.

We will then show that equilibria solutions of \eqref{symPDE} and its reflection are related by an inequality, for small $\epsilon>0$. Then one can extend such inequality for larger $\epsilon$. Lastly one can prove the reversed inequality. The method is better illustrated in the below picture and the following lemmata.
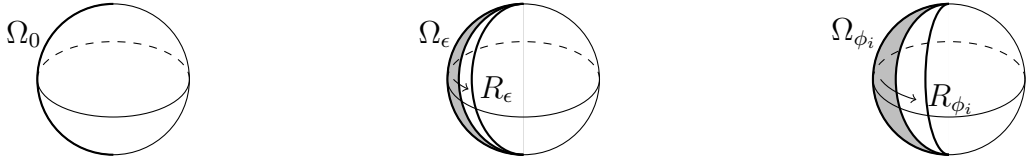
\begin{figure}[!htb]
\minipage{0.3\textwidth}
    \begin{tikzpicture}[scale=1]
    \draw (-1,0) arc (180:360:1cm and 0.5cm);
    \draw[dashed] (-1,0) arc (180:0:1cm and 0.5cm);
    \draw (0,0) circle (1cm);
    \draw[thick] (0,1) arc (90:270:1cm and 1cm);
    \draw (-0.8,0.6) node[left] {$\Omega_0$};
    \end{tikzpicture}
\endminipage\hfill
\minipage{0.3\textwidth}%
    \begin{tikzpicture}[scale=1]
    \fill [fill=lightgray] (0,1) arc (90:270:1cm and 1cm) -- (0,1) arc (90:270:0.7cm and 1cm);
    \fill [fill=white] (0,1) arc (90:270:0.85cm and 1cm) -- (0,1);    
    \draw (-1,0) arc (180:360:1cm and 0.5cm);
    \draw[dashed] (-1,0) arc (180:0:1cm and 0.5cm);
    \draw (0,0) circle (1cm);
    \draw[thick] (0,1) arc (90:270:1cm and 1cm);
    \draw[thick] (0,1) arc (90:270:0.85cm and 1cm);    
    \draw[thick] (0,1) arc (90:270:0.68cm and 1cm);
    \draw (-0.8,0.6) node[left] {$\Omega_{\epsilon}$};
    \draw[->] (-0.92,0) arc (200:222:1cm and 0.4cm) node[right] {$R_{{\epsilon}}$};
    \end{tikzpicture}
\endminipage\hfill
\minipage{0.3\textwidth}
    \begin{tikzpicture}[scale=1]
    \fill [fill=lightgray] (0,1) arc (90:270:1cm and 1cm) -- (0,1) arc (90:270:0.7cm and 1cm);
    \fill [fill=white] (0,1) arc (90:270:0.7cm and 1cm) -- (0,1);
    \draw (-1,0) arc (180:360:1cm and 0.5cm);
    \draw[dashed] (-1,0) arc (180:0:1cm and 0.5cm);
    \draw (0,0) circle (1cm);
    \draw[thick] (0,1) arc (90:270:1cm and 1cm);
    \draw[thick] (0,1) arc (90:270:0.7cm and 1cm);
    \draw (-0.8,0.6) node[left] {$\Omega_{\phi_i}$};
    \draw[thick] (0,1) arc (90:270:0.31cm and 1cm);
    \draw[->] (-0.9,0) arc (200:242:1cm and 0.5cm) node[right] {$R_{\phi_i}$};
    \end{tikzpicture}
\endminipage
\caption{Moving Arc Method: start at $\phi_{i-1}$ and consider sectors $\Omega_{\epsilon}$ nearby it. Reflect along the and $\epsilon$-arc to obtain an inequality. Then extend such inequality for bigger $\epsilon$. Repeat the process starting at $\phi_{i+1}$ and move the arc in the opposite direction.}
\end{figure}



\begin{lem} \textbf{Inequality for Reflections}\label{arbysrot}

Consider $u$ a nonconstant equilibrium of \eqref{symPDE} such that all extrema are axial. If $\phi_i$ is an axial maximum, then
\begin{equation*}
    u(\theta,\phi)\leq u(\theta,R_{\phi_i}(\phi))
\end{equation*}
where $(\theta,\phi)\in[0,\pi] \times [\phi_{i}-\phi^*_i,\phi_{i}]$ and $\phi^*_i:=min \{ |\phi_i-\phi_{i-1}|,|\phi_i-\phi_{i+1}|\}$. 
\end{lem}

Now, we can repeat the process in the above lemma, by moving the arc from $\phi_{i+1}$ in the opposite direction in order to obtain the reversed inequality, yielding the reflection symmetry in the $\phi$ direction. Note such inequalities will only hold for certain subdomains.

\begin{lem} \textbf{Equality for Reflections}\label{arbyseq}

Consider $u$ a non constant equilibrium of \eqref{symPDE} with only axial extrema.  If $\phi_i$ is an axial maximum, then
\begin{equation*}
    u(\theta,\phi)= u(\theta,R_{\phi_i}(\phi))
\end{equation*}
where $(\theta,\phi)\in[0,\pi] \times [\phi_{i}-\phi^*_i,\phi_{i}]$ and $\phi^*_i$ is as before.
\end{lem}

Similarly, the Lemma \ref{arbysrot} holds with reversed inequality, if $\phi_i$ is an axial minimum, and once can repeat the process of Lemma \ref{arbyseq} in order to obtain an equality, and prove the reflectional symmetry along minima.

The main tools to prove such Lemmata are different versions of the maximum principle. Below we state a particular version of the one for thin domains on manifolds, as in Padilla \cite{Padilla}, and the strong unique continuation theorem for manifolds, due to Kazdan \cite{Kazdan88}.
\begin{thm}\emph{\textbf{Maximum Principles}} \label{MPthm}

Consider $w$ a solution of the linear differential equation,
\begin{equation}
    0=a(\theta,\phi)\Delta_{\mathbb{S}^2}w+b(\theta,\phi)w
\end{equation}
where $(\theta,\phi)\in \Omega_\epsilon \subseteq \mathbb{S}^2$ with $\epsilon>0$, and $a,b\in L^\infty_{loc}$ such that $a>0$.
\begin{enumerate}
    \item {\emph{Thin domains:}} Consider $U\subseteq \Omega_\epsilon$ such that the Lebesgue measure $|U|<\mu$ for sufficiently small $\mu$. If $w\leq 0$ on $ \partial U$, then $w\leq 0$ in $U$.
    \item \emph{Strong:} If $w\leq 0$ on $\Omega_\epsilon$, then either $w < 0$ in the interior of $\Omega_\epsilon$ or $w\equiv 0$ in $\overline{\Omega_\epsilon}$.
    \item \emph{Unique continuation:} If $w\equiv 0$ in an open set $U\subseteq \mathbb{S}^2$, then $w\equiv 0$ in $\mathbb{S}^2$.
\end{enumerate}
\end{thm}

In particular, $\Omega_\epsilon$ with $\epsilon>0$ sufficiently small is a thin domain. Also, the sector $\Omega_\epsilon$ is an open set.

In order to prove the strong version, suppose that $w({\theta},{\phi}) =0$ for some $({\theta},{\phi})\in \Omega_\epsilon$. One can calculate the value of a function at a point by its mean value around a geodesic ball $B_r$ in the sphere with geodesic distance $r>0$ small. Hence,
\begin{equation*}
    0=w({\theta},{\phi})=\frac{1}{vol(B_r)}\int_{B_r} w \leq \sup_{B_r}(w)=0
\end{equation*}
and the equality holds if and only if $w\equiv 0$ on $B_r$. One can show that the set $\{ ({\theta},{\phi})\in \Omega_\epsilon \text{ $|$ } w({\theta},{\phi})=0\}$ is nonempty, closed and open in the connected set $\Omega_\epsilon$, hence it must be the full set $\Omega_\epsilon$. See \cite{Evans10} for a detailed exposition of the Euclidean case. 

We proceed to finish the proof of Theorem \ref{symhetthm4} by proving the above lemmata.

\begin{pf} \textbf{of Lemma \ref{arbysrot}.} We firstly prove an inequality for a small reflection, then we extend the inequality for larger reflections.


Recall that the axial extrema are finite, since $f$ is analytic, and so is $u$. See \cite{CaoRammahaTiti00}. Consider three consecutive axial extrema, $\phi_{i-1}<\phi_i<\phi_{i+1}$. As mentioned, we start the moving arc at $\phi_{i-1}$. Suppose that it is an axial minimum, then $\phi_{i}$ is an axial maximum and $\phi_{i+1}$ is an axial minimum. 

Indeed, if $\phi_{i}$ was a also an axial minimum, there would be an axial maximum $\phi_*\in(\phi_{i-1},\phi_{i})$, by the extreme value theorem and noticing all extrema are axial. This contradicts that $\phi_{i-1}$ and $\phi_i$ are consecutive. Similarly for $\phi_{i+1}$.

Therefore,
\begin{equation}\label{pos}
   u(\theta,\phi_{i-1})\leq u(\theta,\phi)
\end{equation} 
for all $\theta\in[0,\pi]$ and $\phi\in (\phi_{i-1},\phi_{i})$. In particular, the inequality must be strict for some value $(\theta,\phi)$, otherwise the function would be locally constant, and due to the unique continuation theorem, globally constant.

Define the difference of the equilibrium solution $u$ of \eqref{symPDE} and itself with reflected angle through
\begin{equation*}
    w^\epsilon(\theta,\phi):=u(\theta,\phi)-u(\theta,R_\epsilon(\phi))    
\end{equation*} 
for $(\theta,\phi)\in\Omega_\epsilon$ with $\epsilon>0$ sufficiently small. Then, $w^\epsilon$ satisfies
\begin{align}
\begin{cases}\label{ellipOP}
    0=Lw^\epsilon+b(\theta,\phi)w^\epsilon &\text{ on } \Omega_\epsilon\\
    w^\epsilon\leq 0  &\text{ on } \partial\Omega_\epsilon
\end{cases}
\end{align}
where $L:=a(\theta,\phi)\Delta_{\mathbb{S}^2}$ is an elliptic operator defined below.

Indeed, the boundary values are $w^\epsilon(\theta,\epsilon)=0$ due to the reflection, and $w^\epsilon(\theta,0)\leq 0$ since $\phi_{i-1}=0$ is an axial minimum, as in \eqref{pos}. Moreover, the fundamental theorem of calculus implies that the difference of two solutions of \eqref{symPDE} satisfy
\begin{equation*}
    0= \int^1_0 \frac{d}{d\tau} \left[ F(u^\tau,\Delta_{\mathbb{S}^2} u^\tau) \right] d\tau
\end{equation*}
with $u^\tau:=\tau u(\theta,\phi)+(1-\tau)u(\theta,R_\epsilon (\phi))$. 
Calculating the $\tau$-derivative,
\begin{equation*}
    0 = a(\theta,\phi)\Delta_{\mathbb{S}^2}w^\epsilon+b(\theta,\phi)w^\epsilon
\end{equation*}
where the bounded coefficients are given by the Hadamard formulas
\begin{align*}
    a(\theta,\phi)&:=\int_0^1 F_q(u^\tau,\Delta_{\mathbb{S}^2}u^\tau))d\tau \\
    b(\theta,\phi)&:=\int_0^1 F_u(u^\tau,\Delta_{\mathbb{S}^2}u^\tau))d\tau d\tau 
\end{align*}
where $(u,q):=(u,\Delta_{\mathbb{S}^2}u)$.

The maximum principle for thin domains in Theorem \ref{MPthm} implies that
\begin{equation}\label{locARG}
    w^\epsilon\leq 0
\end{equation}
in $\Omega_\epsilon$ for $\epsilon>0$ sufficiently small. 



This inequality is now extended for larger $\epsilon$. Let
\begin{equation*}
    \epsilon^*:=\sup\{ \epsilon\in [0,2\pi] \text{ $|$ } w^\epsilon\leq 0  \text{ in } \Omega_{\epsilon} \text{ for all } \epsilon\leq \epsilon^*\}. 
\end{equation*}

Note $\epsilon^*>0$, using the maximum principle for thin domains as in \eqref{locARG}. Also, $\epsilon^*\leq \phi_{i}$, since this is an axial maximum and hence $w^{\phi_{i}+\delta}(\theta,\phi)> 0$ for some $(\theta,\phi)\in\Omega_{\phi_{i}+\delta} \backslash \Omega_{\phi_{i}}$ with $\delta>0$ sufficiently small, by applying the same arguments as in \eqref{pos} for a maximum. 

We claim that $\epsilon^*=\phi_{i}$. Suppose that $\epsilon^*<\phi_{i}$. In order to contradict its maximality, it is shown that
\begin{equation}\label{RandomIneq}
    w^{\epsilon^*+\delta}\leq 0    
\end{equation} 
in $\Omega_{\epsilon^*+\delta}$ for some $\delta>0$. For such extension argument, we adapt the method using the unique continuation theorem, as Pol{\'a}\v{c}ik \cite{Polacik12}.

By hypothesis, we have that $w^{\epsilon^*}\leq 0$ in $\Omega_{\epsilon^*}$. Moreover, the strong maximum principle in Theorem \ref{MPthm} guarantees 
\begin{equation*}
    w^{\epsilon^*}<0    
\end{equation*}
in the interior of $\Omega_{\epsilon^*}$. The other possibility, that $w^{\epsilon^*}\equiv 0$ in $\overline{\Omega_{\epsilon^*}}$ would imply that $w^{\epsilon^*}\equiv 0$ in $\mathbb{S}^2$ due to the unique continuation theorem. Hence $u$ would be a constant, which is a contradiction. 

In particular, 
\begin{equation*}
    w^{\epsilon^*}<0    
\end{equation*}
for a compact subset 
$K$ of ${\Omega_{\epsilon^*}}$ such that $|\Omega_{\epsilon^*}\backslash K | < \frac{\mu}{2}$. Due to continuity in the parameter and that $K$ is compact,
\begin{equation}\label{wtf1000}
    w^{\epsilon^*+\delta}< 0    
\end{equation}
in $K$.

To prove \eqref{RandomIneq}, and contradict the maximality of $\epsilon^*$, it is enough to prove an inequality for the remaining set $\Omega_{\epsilon^*+\delta}\backslash K$. In this set, $w^{\epsilon^*+\delta}$ satisfy
\begin{align}\label{enoughofthisbullshit}
\begin{cases}
    0=Lw^{\epsilon^*+\delta}+c(\theta,\phi)w^{\epsilon^*+\delta} &\text{ on } \Omega_{\epsilon^*+\delta}\backslash K\\
    w^{\epsilon^*+\delta}\leq 0  &\text{ on } \partial (\Omega_{\epsilon^*+\delta}\backslash K)=\partial \Omega_{\epsilon^*+\delta}\cup \partial K
\end{cases}
\end{align}
with $L$ and $c(\theta,\phi)$ as in \eqref{ellipOP}. Note that one can choose $\delta>0$ small such that $|\Omega_{\epsilon^*+\delta}\backslash \Omega_{\epsilon^*}| \leq \frac{\mu}{2}$, and hence $|\Omega_{\epsilon^*+\delta}\backslash K| \leq \mu$. Also, since $K$ is a compact subset of the open set $\Omega_{\epsilon^*+\delta}$, then the boundary $\partial(\Omega_{\epsilon^*+\delta}\backslash K)=\partial\Omega_{\epsilon^*+\delta}\cup\partial K$. Moreover, the boundary values of \eqref{enoughofthisbullshit} in $\partial K$ are obtained because of the inequality \eqref{wtf1000}, $w^{\epsilon^*+\delta}(\theta,\epsilon^*+\delta)=0$ due to the reflection at the $(\epsilon^*+\delta)$-line, and $w^{\epsilon^*+\delta}(\theta,0)\leq 0$ since the solution only has leveled axial extrema, where $(\theta,0)$ is an axial minimum, and $\phi_{i}$ is closer to $\phi_{i-1}$.

Since $\Omega_{\epsilon^*+\delta}\backslash K$ is still a thin domain for $\delta>0$ small, the maximum principle implies that
\begin{equation}\label{yetanotherone}
    w^{\epsilon^*+\delta}\leq 0    
\end{equation}
in $\Omega_{\epsilon^*+\delta}\backslash K$. 

Combining the inequality \eqref{yetanotherone} in $\Omega_{\epsilon^*+\delta}\backslash K$ with the inequality \eqref{wtf1000} in $K$, yields the desired inequality \eqref{RandomIneq} in the whole set $\Omega_{\epsilon^*+\delta}$. This contradicts the maximality of $\epsilon^*$, yielding
\begin{equation*}
    w^{\phi_i}\leq 0
\end{equation*} 
in $\Omega_{\phi_i}$.

\begin{flushright}
	$\blacksquare$
\end{flushright}
\end{pf}

\begin{pf} \textbf{of Lemma \ref{arbyseq}.} 

Consider three consecutive leveled axial extrema $\phi_{i-1}<\phi_{i}<\phi_{i+1}$ such that $\phi_i$ is a maximum and the other two are minima. From the Lemma \ref{arbysrot},
\begin{equation}\label{HJ}
    u(\theta,{\phi})\leq u(\theta,R_{{\phi}_i}({\phi}))
\end{equation}
with $(\theta,{\phi})\in[0,\pi] \times [{\phi}_{i-1},{\phi}_i]$.

It is proven the reversed inequality. This is done by moving the arc on the reversed orientation of $\phi$, starting at the minimum $\phi_{i+1}$. Indeed, consider the change of variable $\tilde{\phi}(\phi):=\phi_{i+1}-\phi$ with starting axis $\tilde{\phi}(\phi_{i+1})=0$.

The analogous of the condition \eqref{pos} is that $\tilde{\phi}=0$ is an axial minimum, namely
\begin{equation*}
    u(\theta,0)\leq u(\theta,\tilde{\phi})    
\end{equation*}
for $\tilde{\phi}\in (0,\tilde{\phi}_i)$ and $\tilde{\phi}_i:=\tilde{\phi}(\phi_i)=\phi_{i+1}-\phi_i$ is where the maximum $\phi_i$ lies in the $\tilde{\phi}$ coordinates. 

Hence, one can apply the Lemma \ref{arbysrot} in the new variable $\tilde{\phi}$,
\begin{equation}\label{HJ2}
    u(\theta,\tilde{\phi})\leq u(\theta,R_{\tilde{\phi}_i}(\tilde{\phi}))
\end{equation}
with $(\theta,\tilde{\phi})\in[0,2\pi] \times [0,\tilde{\phi}_i]$.

We want to compare \eqref{HJ} and \eqref{HJ2}, but they are valid in different domains. Consider $\phi\in [\phi_{i-1},\phi_{i}]$, yielding \eqref{HJ}. Then, either $R_{{\phi}_i}({\phi})$ is in $  [{\phi}_i,{\phi}_{i+1}]$ or $[{\phi}_{i-1},{\phi}_{i}]$, depending which interval is bigger. In the former case, one can use \eqref{HJ2} with $\tilde{\phi}=R_{{\phi}_i}({\phi})$, whereas in the latter, one can use \eqref{HJ} again with $R_{{\phi}_i}({\phi})$ instead of $\phi$. Both cases yield
\begin{equation*}
    u(\theta,\phi)\leq u(\theta,R_{\phi_i}(\phi))\leq u(\theta,R_{\tilde{\phi}_i}(R_{\phi_i}(\phi)))
\end{equation*}
with $\phi\in [\phi_{i-1},\phi_{i}]$. 

Since $\phi_i$ and $\tilde{\phi}_i$ denote the same point in different coordinates, the composition of reflections $R_{\tilde{\phi}_i}(R_{\phi_i}(\phi))=\phi$. Hence, 
\begin{equation*}
    u(\theta,\phi)= u(\theta,R_{\phi_i}(\phi))
\end{equation*}
for $\phi\in [\phi_{i-1},\phi_{i}]$. 

In particular,
\begin{equation*}
    u(\theta,\phi_{i-1})= u(\theta,R_{\phi_i}(\phi_{i-1})).
\end{equation*}

Then $R_{\phi_i}(\phi_{i-1})$ is another leveled and axial minimum, since all extrema are leveled and axial. Moreover, we chose three consecutive extrema $\phi_{i-1}<\phi_{i}<\phi_{i+1}$, hence $R_{\phi_i}(\phi_{i-1})=\phi_{i+1}$ and $\phi_i$ is the midpoint between $\phi_{i-1}$ and $\phi_{i+1}$.
\begin{flushright}
	$\blacksquare$
\end{flushright}
\end{pf}

\section{Discussion}\label{sec:conclusion}

We comment on possible generalizations and extensions of Theorem \ref{symhetthm4}. 

In particular, the hypothesis can be weakened yielding similar results. If the leveled assumption is dropped out, one might be able to recover a result for a subdomain in $\phi$. Or if the extrema are not axial, but are curves from the north pole to the south pole, then a similar result should hold, using the methods of Pol{\'a}\v{c}ik \cite{Polacik12}. 
Also, if the regularity of $F\in C^1$, it should be possible to obtain the same result, since analyticity was only used to obtain finitely many axial extrema.

Note that we do not consider $F$ depending on the gradient term $\nabla_{\mathbb{S}^2}$ with the Levi-Civita connection, since $u$ and its reflection with respect to an extrema would satisfy different equations. We can consider $F(u,\nabla_{\mathbb{S}^2}u,\Delta_{\mathbb{S}^2}u)$ by considering additional hypothesis on $F$.

Another possible extension of the main theorem is a symmetry property for parabolic equations. Indeed, similar results for positive solutions of parabolic equations on the ball with Dirichlet boundary were obtained by Pol{\'a}\v{c}ik and Hess \cite{HessPolacik94}, and Babin \cite{Babin94}. In particular, it is proven that the global attractor $\mathcal{A}$ restricted to the subspace of positive solutions consists of radially symmetric solutions. Hence, the result for elliptic equations can be seen as a particular case of the parabolic version, since equilibria are in the attractor.

Note the moving arc method here presented is not limited to the sphere. Depending on the Laplace-Beltrami, the results here can be adapted to other domains. For general manifolds, one needs an assumption on convexity in the candidate direction for symmetry, such as geodesic convexity. Further investigation is being carried for the torus and the hyperbolic disk.

For instance, Theorem \ref{symhetthm4} also holds for reflections of the angle $\phi_{1}$ of the sphere $\mathbb{S}^n$ with $n\geq 1$. Indeed, in geodesic coordinates $(\theta,\phi) \in {\mathbb{S}^n}$, where $\theta\in [0,\pi]$ is the geodesic radius starting from the north pole, and $\phi=(\phi_1,...,\phi_{n-1})\in \mathbb{S}^{n-1}$ is the angle coordinate from a particular distance of the north pole. 

In particular, for $n=1$, the axial extrema hypothesis is superfluous and the only remaining condition is that extrema are leveled. This implies the symmetry of the equilibria with leveled extrema in the attractor constructed by Fiedler, Rocha and Wolfrum \cite{FiedlerRochaWolfrum04}.

\end{document}